\documentclass{amsart}
\usepackage{amsmath,amssymb}
\usepackage{yhmath}

\newcommand{\bdis}{\begin{displaymath}}
\newcommand{\edis}{\end{displaymath}}
\newcommand{\be}{\begin{equation}}
\newcommand{\ee}{\end{equation}}
\newcommand{\mbb}{\mathbb}
\newcommand{\mcal}{\mathcal}

\newcommand{\vp}{\varphi}

\newcommand{\zf}{\zeta\left(\frac{1}{2}+it\right)}

\newcommand{\FR}{\frac{x^n+y^n}{z^n}}

\theoremstyle{definition}

\theoremstyle{remark}
\newtheorem{remark}[]{Remark}

\newtheorem*{mydef11}{{\bf Theorem 1}}

\newtheorem*{mydef12}{{\bf Theorem 2}}

\newtheorem*{mydef13}{{\bf Theorem 3}}

\newtheorem*{mydef14}{{\bf Theorem 4}}

\newtheorem*{mydef51}{{\bf Lemma 1}}

\newtheorem*{mydef52}{{\bf Lemma 2}}

\newtheorem*{mydef53}{{\bf Lemma 3}}

\newtheorem*{mydef54}{{\bf Lemma 4}}

\newtheorem*{mydef55}{{\bf Lemma 5}}

\newtheorem*{mydef81}{{\bf Property 1}}

\newtheorem*{mydef82}{{\bf Property 2}}

\numberwithin{equation}{section}

\begin{document}

\title[Jacob's ladders and new equivalents \dots]{Jacob's ladders and new equivalents of the Fermat-Wiles theorem connected with some cross-bred of the formulae of Hardy-Littlewood-Ingham (1926) and of Ingham (1926)}

\author{Jan Moser}

\address{Department of Mathematical Analysis and Numerical Mathematics, Comenius University, Mlynska Dolina M105, 842 48 Bratislava, SLOVAKIA}

\email{jan.mozer@fmph.uniba.sk}


\begin{abstract}
The main result of this paper is new formula connecting certain $zeta$-integral on the critical line with a $\zeta$-integral in the critical strip. Further, a kind of cross-breeding of the Hardy-Littlewood-Ingham formula and Ingham formula produces new $\zeta$-equivalent of the Fermat-Wiles theorem. 
\end{abstract}
\maketitle

\section{Introduction} 

\subsection{} 

Let us remind the two sets of increments 
\be \label{1.1} 
\begin{split}
& \left\{ \int_{\overset{r-1}{T}(T)}^{\overset{r}{T}(T)}\left|\zf\right|^2{\rm d}t\right\},\ 
\left\{ \int_{\overset{r-1}{T}(T)}^{\overset{r}{T}(T)}\left|\zeta(\sigma+it)\right|^2{\rm d}t\right\}, \\ 
& r=1,\dots, k
\end{split} 
\ee 
for every fixed $k\in\mbb{N}$ and every fixed, sufficiently big $T>0$, of the Hardy-Littlewood integrals: (1918) 
\be \label{1.2} 
J(T)=\int_0^T\left|\zf\right|^2{\rm d}t,\ T\to\infty, 
\ee  
and (1922) 
\be \label{1.3} 
J_1(T,\sigma)=\int_1^T|\zeta(\sigma+it)|^2{\rm d}t,\ \sigma\geq \frac 12+\epsilon, 
\ee  
respectively, for every fixed $\epsilon>0$. 

In this paper, that is a continuation of the paper's series \cite{8} -- \cite{15}, we prove firstly the following new result 
\be \label{1.4} 
\frac{ \int_{\overset{r-1}{T}(T)}^{\overset{r}{T}(T)}\left|\zf\right|^2{\rm d}t}{\int_{\overset{r-1}{T}(T)}^{\overset{r}{T}(T)}\left|\zeta(\sigma+it)\right|^2{\rm d}t}\sim \frac{1}{\zeta(2\sigma)}\ln\overset{r-1}{T},\ T\to\infty 
\ee 
for every fixed: 
\bdis 
k\in\mbb{N},\ \sigma\geq \frac 12+\epsilon. 
\edis 

\begin{remark}
Certain non-local interaction of two integrals 
\bdis 
\int_{\overset{r-1}{T}(T)}^{\overset{r}{T}(T)}\left|\zf\right|^2{\rm d}t,\ \int_{\overset{r-1}{T}(T)}^{\overset{r}{T}(T)}\left|\zeta(\sigma+it)\right|^2{\rm d}t
\edis 
is expressed by the formula (\ref{1.4}) in the following sense 
\bdis 
t\in [\overset{r-1}{T}(T),\overset{r}{T}(T)] \ \Rightarrow \ \left\{\frac 12+it\right\}\cap \{\sigma+it\}=\emptyset, 
\edis 
that is the segments of integrations in (\ref{1.4}) constitute a disconnected set, or sets $r=1,\dots,k$. 
\end{remark}  

\subsection{} 

Let us remind the classical Hardy-Littlewood-Ingham formula 
\be \label{1.5} 
\int_0^T\left|\zf\right|^2{\rm d}t=T\ln T-(1+\ln2\pi-2c)T+\mcal{O}(T^{1/2}\ln T) 
\ee 
(for our purpose, this expression for the error term is quite enough). 

Next, we give here some variant of \emph{cross-breeding} of the formula (\ref{1.5}) with the formula (\ref{1.4}) to obtain the following result: It is the cross-breed functional of the Hardy-Littlewood-Ingham formula 
\be \label{1.6} 
\begin{split}
& \lim_{\tau\to\infty}\frac{1}{\tau}\int_1^{\frac{x}{\zeta(2\sigma)}\tau}\left|\zf\right|^2{\rm d}t \times \\ 
& \int_{\frac{x}{\zeta(2\sigma)}\tau}^{[\frac{x}{\zeta(2\sigma)}\tau]^1}|\zeta(\sigma+it)|^2{\rm d}t \times \\ 
& \left\{\int_{\frac{x}{\zeta(2\sigma)}\tau}^{[\frac{x}{\zeta(2\sigma)}\tau]^1}\left|\zf\right|^2{\rm d}t\right\}^{-1} = x 
\end{split}
\ee  
for every fixed 
\bdis 
x>0,\ \sigma\geq\frac 12+\epsilon, 
\edis  
and, of course, 
\bdis 
[Y]^1=\vp_1^{_1}(Y). 
\edis 

\begin{remark}
It is then clear, that in the special case 
\be \label{1.7} 
x\to \frac{x^n+y^n}{z^n},\ x,y,z,n\in\mbb{N},\ n\geq 3, 
\ee  
the new $\zeta$-equivalent of the Fermat-Wiles theorem follows from (\ref{1.6}) for every fixed $\sigma\geq \frac 12+\epsilon$. 
\end{remark} 

\subsection{} 

Next, let us remind the classical Ingham formula, see \cite{2}, 
\be \label{1.8} 
\int_1^T\left|\zf\right|^4{\rm d}t=\frac{1}{2\pi^2}T\ln^4T+\mcal{O}(T\ln^3T),\ T\to\infty. 
\ee 
The cross-breeding of this formula with the formula (\ref{1.4}) gives the following: the $\zeta$-condition 
\be \label{1.9} 
\begin{split}
& \lim_{\tau\to\infty}\frac{1}{\tau}\int_1^{\frac{2\pi^2}{\zeta^4(2\sigma)}\FR\tau}\left|\zf\right|^4{\rm d}t \times \\ 
& \left\{\int_{\frac{2\pi^2}{\zeta^4(2\sigma)}\FR\tau}^{[\frac{2\pi^2}{\zeta^4(2\sigma)}\FR\tau]^1}|\zeta(\sigma+it)|^2\right\}^4 \times \\ 
& \left\{\int_{\frac{2\pi^2}{\zeta^4(2\sigma)}\FR\tau}^{[\frac{2\pi^2}{\zeta^4(2\sigma)}\FR\tau]^1}\left|\zf\right|^2\right\}^{-4} \not= 1
\end{split}
\ee 
on the class of all Fermat's rationals, see (\ref{1.7}), represents the next $\zeta$-equivalent of the Fermat-Wiles theorem for every fixed $\sigma\geq\frac 12+\epsilon$.

\section{Jacob's ladders: notions and basic geometrical properties}  

\subsection{} 

In this paper we use the following notions of our works \cite{3} -- \cite{7}: 
\begin{itemize}
\item[{\tt (a)}] Jacob's ladder $\vp_1(T)$, 
\item[{\tt (b)}] direct iterations of Jacob's ladders 
\bdis 
\begin{split}
	& \vp_1^0(t)=t,\ \vp_1^1(t)=\vp_1(t),\ \vp_1^2(t)=\vp_1(\vp_1(t)),\dots , \\ 
	& \vp_1^k(t)=\vp_1(\vp_1^{k-1}(t))
\end{split}
\edis 
for every fixed natural number $k$, 
\item[{\tt (c)}] reverse iterations of Jacob's ladders 
\be \label{2.1}  
\begin{split}
	& \vp_1^{-1}(T)=\overset{1}{T},\ \vp_1^{-2}(T)=\vp_1^{-1}(\overset{1}{T})=\overset{2}{T},\dots, \\ 
	& \vp_1^{-r}(T)=\vp_1^{-1}(\overset{r-1}{T})=\overset{r}{T},\ r=1,\dots,k, 
\end{split} 
\ee   
where, for example, 
\be \label{2.2} 
\vp_1(\overset{r}{T})=\overset{r-1}{T}
\ee  
for every fixed $k\in\mbb{N}$ and every sufficiently big $T>0$. We also use the properties of the reverse iterations listed below.  
\be \label{2.3}
\overset{r}{T}-\overset{r-1}{T}\sim(1-c)\pi(\overset{r}{T});\ \pi(\overset{r}{T})\sim\frac{\overset{r}{T}}{\ln \overset{r}{T}},\ r=1,\dots,k,\ T\to\infty,  
\ee 
\be \label{2.4} 
\overset{0}{T}=T<\overset{1}{T}(T)<\overset{2}{T}(T)<\dots<\overset{k}{T}(T), 
\ee 
and 
\be \label{2.5} 
T\sim \overset{1}{T}\sim \overset{2}{T}\sim \dots\sim \overset{k}{T},\ T\to\infty.   
\ee  
\end{itemize} 

\begin{remark}
	The asymptotic behaviour of the points 
	\bdis 
	\{T,\overset{1}{T},\dots,\overset{k}{T}\}
	\edis  
	is as follows: at $T\to\infty$ these points recede unboundedly each from other and all together are receding to infinity. Hence, the set of these points behaves at $T\to\infty$ as one-dimensional Friedmann-Hubble expanding Universe. 
\end{remark}  

\subsection{} 

Let us remind that we have proved\footnote{See \cite{8}, (3.4).} the existence of almost linear increments 
\be \label{2.6} 
\begin{split}
& \int_{\overset{r-1}{T}}^{\overset{r}{T}}\left|\zf\right|^2{\rm d}t\sim (1-c)\overset{r-1}{T}, \\ 
& r=1,\dots,k,\ T\to\infty,\ \overset{r}{T}=\overset{r}{T}(T)=\vp_1^{-r}(T)
\end{split} 
\ee 
for the Hardy-Littlewood integral (1918) 
\be \label{2.7} 
J(T)=\int_0^T\left|\zf\right|^2{\rm d}t. 
\ee  

For completeness, we give here some basic geometrical properties related to Jacob's ladders. These are generated by the sequence 
\be \label{2.8} 
T\to \left\{\overset{r}{T}(T)\right\}_{r=1}^k
\ee 
of reverse iterations of of the Jacob's ladders for every sufficiently big $T>0$ and every fixed $k\in\mbb{N}$. 

\begin{mydef81}
The sequence (\ref{2.8}) defines a partition of the segment $[T,\overset{k}{T}]$ as follows 
\be \label{2.9} 
|[T,\overset{k}{T}]|=\sum_{r=1}^k|[\overset{r-1}{T},\overset{r}{T}]|
\ee 
on the asymptotically equidistant parts 
\be \label{2.10} 
\begin{split}
& \overset{r}{T}-\overset{r-1}{T}\sim \overset{r+1}{T}-\overset{r}{T}, \\ 
& r=1,\dots,k-1,\ T\to\infty. 
\end{split}
\ee 
\end{mydef81} 

\begin{mydef82}
Simultaneously with the Property 1, the sequence (\ref{2.8}) defines the partition of the integral 
\be \label{2.11} 
\int_T^{\overset{k}{T}}\left|\zf\right|^2{\rm d}t
\ee 
into the parts 
\be \label{2.12} 
\int_T^{\overset{k}{T}}\left|\zf\right|^2{\rm d}t=\sum_{r=1}^k\int_{\overset{r-1}{T}}^{\overset{r}{T}}\left|\zf\right|^2{\rm d}t, 
\ee 
that are asymptotically equal 
\be \label{2.13} 
\int_{\overset{r-1}{T}}^{\overset{r}{T}}\left|\zf\right|^2{\rm d}t\sim \int_{\overset{r}{T}}^{\overset{r+1}{T}}\left|\zf\right|^2{\rm d}t,\ T\to\infty. 
\ee 
\end{mydef82} 

It is clear, that (\ref{2.10}) follows from (\ref{2.3}) and (\ref{2.5}) since 
\be \label{2.14} 
\overset{r}{T}-\overset{r-1}{T}\sim (1-c)\frac{\overset{r}{T}}{\ln \overset{r}{T}}\sim (1-c)\frac{T}{\ln T},\ r=1,\dots,k, 
\ee  
while our eq. (\ref{2.13}) follows from (\ref{2.6}) and (\ref{2.5}). 

\section{Theorem about non-local interaction of the two integrals} 

\subsection{} 

Let us remind the following Hardy-Littlewood formula (1922) 
\be \label{3.1} 
\int_1^T|\zeta(\sigma+it)|^2{\rm d}t=\zeta(2\sigma)T+\mcal{O}(T^{1-\epsilon}\ln T),\ T\to\infty, 
\ee  
for every fixed: 
\be \label{3.2} 
\epsilon>0,\ \sigma\geq \frac 12+\epsilon. 
\ee 

\begin{remark}
Hardy and Littlewood have proved just non-trivial part of this formula, namely the part that is true for the critical strip, comp. \cite{16}, pp. 30, 31. 
\end{remark} 

Since 
\be \label{3.3} 
\overset{r}{T}-\overset{r-1}{T}\sim (1-c)\frac{T}{\ln T},\ r=1,\dots,k,\ T\to\infty, 
\ee  
and, comp. (\ref{3.1}), 
\be \label{3.4} 
T^{1-\epsilon}\ln T=o\left(\frac{T}{\ln T}\right),\ T\to\infty, 
\ee  
it is true the following. 

\begin{mydef51}
The formula 
\be \label{3.5} 
\begin{split}
& \int_{\overset{r-1}{T}}^{\overset{r}{T}}|\zeta(\sigma+it)|^2{\rm d}t=\zeta(2\sigma)(\overset{r}{T}-\overset{r-1}{T})+\mcal{O}(T^{1-\epsilon}\ln T), \\ 
& r=1,\dots,k,\ T\to\infty 
\end{split} 
\ee 
is the asymptotic one for every fixed $\sigma=\frac{1}{2}+\epsilon$ and $k\in\mbb{N}$. 
\end{mydef51} 

We use in this section also our almost linear formula\footnote{See \cite{8}, (3.4).} 
\be \label{3.6} 
\int_{\overset{r-1}{T}}^{\overset{r}{T}}\left|\zf\right|^2{\rm d}t=(1-c)\overset{r-1}{T}+\mcal{O}(T^{\frac 13+\delta}) 
\ee  
for every fixed $\delta>0$ and $k\in\mbb{N}$. 

\subsection{} 

Now, we have\footnote{See (\ref{3.5}), (\ref{3.6}) and also (\ref{2.5}).} 
\be \label{3.7} 
\frac{\int_{\overset{r-1}{T}}^{\overset{r}{T}}\left|\zf\right|^2{\rm d}t}{\int_{\overset{r-1}{T}}^{\overset{r}{T}}|\zeta(\sigma+it)|^2{\rm d}t}=
\frac{(1-c)\overset{r-1}{T}}{\zeta(2\sigma)(\overset{r}{T}-\overset{r-1}{T})}
\frac{1+\mcal{O}(T^{-\frac 23+\delta})}{1+\mcal{O}(T^{-\epsilon}\ln^2T)}\overset{\star}{=}
\ee 
Next, we use our more exact formula for $\overset{r}{T}-\overset{r-1}{T}$\footnote{Comp. \cite{3}, p. 420.} 
\bdis 
\overset{r}{T}-\overset{r-1}{T}=(1-c)\frac{\overset{r}{T}}{\ln\overset{r}{T}}\left\{1+\mcal{O}\left(\frac{1}{\ln T}\right)\right\}, 
\edis 
(after the substitutions\footnote{Comp. (\ref{2.2}) and (\ref{2.5}).} 
\bdis 
\frac{1}{2}\vp(T)=\vp_1(T),\ T\to \overset{r}{T},\ \vp_1(\overset{r}{T})=\overset{r-1}{T}
\edis 
in the original formula in \cite{3}) in (\ref{3.7}), that gives 
\be\label{3.8} 
\begin{split}
& \overset{\star}{=}\frac{1}{\zeta(2\sigma)}\frac{\overset{r-1}{T}}{\overset{r}{T}}\ln \overset{r}{T}\left\{1+\mcal{O}\left(\frac{1}{\ln T}\right)\right\}\left\{1+\mcal{O}(T^{-\epsilon}\ln^2T)\right\}= \\ 
& = \frac{1}{\zeta(2\sigma)}\frac{\overset{r-1}{T}}{\overset{r}{T}}\ln \overset{r}{T}\left\{1+\mcal{O}\left(\frac{1}{\ln T}\right)\right\}\overset{\Diamond}{=} 
\end{split}
\ee  
but\footnote{See sec. 2.} 
\be \label{3.9} 
\begin{split}
& \frac{\overset{r-1}{T}}{\overset{r}{T}}=1-\frac{\overset{r}{T}-\overset{r-1}{T}}{\overset{r}{T}}=1+\mcal{O}\left(\frac{1}{\ln T}\right), \\ 
& \ln\overset{r}{T}=\ln\overset{r-1}{T}+\ln\frac{\overset{r}{T}}{\overset{r-1}{T}}=\ln \overset{r-1}{T}+\mcal{O}\left(\frac{1}{\ln T}\right), 
\end{split}
\ee 
that gives finally 
\be \label{3.10} 
\overset{\Diamond}{=}\frac{1}{\zeta(2\sigma)}\ln\overset{r-1}{T}+\mcal{O}(1). 
\ee  

Consequently, we have the following result. 

\begin{mydef11}
It is true, that 
\be\label{3.11} 
\frac{\int_{\overset{r-1}{T}}^{\overset{r}{T}}\left|\zf\right|^2{\rm d}t}{\int_{\overset{r-1}{T}}^{\overset{r}{T}}|\zeta(\sigma+it)|^2{\rm d}t}=\frac{1}{\zeta(2\sigma)}\ln\overset{r-1}{T}+\mcal{O}(1),\ r=1,\dots, k,\ T\to\infty 
\ee 
for every fixed $k\in\mbb{N}$ and $\sigma=\frac 12+\epsilon$. 
\end{mydef11} 

\section{Some cross-breeds of the Hardy-Littlewood-Ingham formula and the next equivalent of the Fermat-Wiles theorem} 

\subsection{} 

We choose the classical Hardy-Littlewood-Ingham formula, see \cite{2}, comp. \cite{8}, (2.1) -- (2.7), 
\be \label{4.1} 
\int_1^T\left|\zf\right|^2{\rm d}t=T\ln T-(1+\ln2\pi-2c)T+\mcal{O}(T^{\frac 12}\ln T)
\ee  
as the basic formula for this section. We give here a variant of cross-breeding of the formula (\ref{4.1}) and the formula (\ref{3.11}). 

To attain our goal we use the formula (\ref{4.1}) in the following form 
\be \label{4.2} 
\int_1^T\left|\zf\right|^2{\rm d}t=T\ln T\left\{1+\mcal{O}\left(\frac{1}{\ln T}\right)\right\},\ T\to\infty. 
\ee  
Since 
\bdis 
\frac{1}{\zeta(2\sigma)}\ln T + \mcal{O}(1)=\frac{\ln T}{\zeta(2\sigma)}\left\{1+\mcal{O}\left(\frac{1}{\ln T}\right)\right\}, 
\edis  
then we have (see (\ref{3.11}), $r=1$) 
\be \label{4.3} 
\begin{split}
& T\ln T=\zeta(2\sigma)T\frac{\ln T}{\zeta(2\sigma)}= \\ 
& \zeta(2\sigma)T\frac{\int_T^{\overset{1}{T}}|\zf|^2{\rm d}t}{\int_T^{\overset{1}{T}}|\zeta(\sigma+it)|^2{\rm d}t}\left\{1+\mcal{O}\left(\frac{1}{\ln T}\right)\right\},
\end{split}
\ee 
that gives us the following 
\be \label{4.4} 
\begin{split}
& \int_1^T\left|\zf\right|^2{\rm d}t\times \frac{\int_T^{\overset{1}{T}}|\zeta(\sigma+it)|^2{\rm d}t}{\int_T^{\overset{1}{T}}|\zf|^2{\rm d}t}= \\ 
& \zeta(2\sigma)T\left\{1+\mcal{O}\left(\frac{1}{\ln T}\right)\right\}. 
\end{split}
\ee  

\subsection{} 

Now, the substitution 
\be \label{4.5} 
T=\frac{x}{\zeta(2\sigma)}\tau,\ \{T\to\infty\} \ \Leftrightarrow \ \{\tau\to\infty\}
\ee 
in the eq. (\ref{4.4}) gives the following functional (as the cross-bred  of the Hardy-Littlewood-Ingham formula). 

\begin{mydef52}
\be \label{4.6} 
\begin{split}
& \lim_{\tau\to\infty}\frac{1}{\tau}\int_1^{\frac{x}{\zeta(2\sigma)}\tau}\left|\zf\right|^2{\rm d}t\times \\ 
& \int_{\frac{x}{\zeta(2\sigma)}\tau}^{[\frac{x}{\zeta(2\sigma)}\tau]^1}|\zeta(\sigma+it)|^2{\rm d}t\times \\ 
& \left\{ 
\int_{\frac{x}{\zeta(2\sigma)}\tau}^{[\frac{x}{\zeta(2\sigma)}\tau]^1}\left|\zf\right|^2{\rm d}t
\right\}^{-1}=x 
\end{split}
\ee  
for every fixed 
\bdis 
x>0, \ \sigma\geq \frac 12+\epsilon. 
\edis 
\end{mydef52} 

Next, in the special case\footnote{See (\ref{1.7}).} 
\bdis 
x\to \FR
\edis  
we obtain the following Lemma from (\ref{4.6}). 

\begin{mydef53}
\be \label{4.7} 
\begin{split}
& \lim_{\tau\to\infty}\frac{1}{\tau}\int_1^{\FR\frac{\tau}{\zeta(2\sigma)}}\left|\zf\right|^2{\rm d}t\times \\ 
& \int_{\FR\frac{\tau}{\zeta(2\sigma)}}^{[\FR\frac{\tau}{\zeta(2\sigma)}]^1}|\zeta(\sigma+it)|^2{\rm d}t\times \\ 
& \left\{\int_{\FR\frac{\tau}{\zeta(2\sigma)}}^{[\FR\frac{\tau}{\zeta(2\sigma)}]^1}\left|\zf\right|^2{\rm d}t\right\}^{-1}=\FR
\end{split}
\ee 
for every fixed Fermat's rational and every fixed 
\bdis 
\sigma\geq \frac 12+\epsilon.
\edis 
\end{mydef53} 

Consequently, we have the following result. 

\begin{mydef12}
The $\zeta$-condition 
\be \label{4.8} 
\begin{split}
& \lim_{\tau\to\infty}\frac{1}{\tau}\int_1^{\FR\frac{\tau}{\zeta(2\sigma)}}\left|\zf\right|^2{\rm d}t\times \\ 
& \int_{\FR\frac{\tau}{\zeta(2\sigma)}}^{[\FR\frac{\tau}{\zeta(2\sigma)}]^1}|\zeta(\sigma+it)|^2{\rm d}t\times \\ 
& \left\{\int_{\FR\frac{\tau}{\zeta(2\sigma)}}^{[\FR\frac{\tau}{\zeta(2\sigma)}]^1}\left|\zf\right|^2{\rm d}t\right\}^{-1}\not= 1
\end{split}
\ee 
on the class of all Fermat's rationals represents the next $\zeta$-equivalent of the Fermat-Wiles theorem for every fixed $\sigma\geq\frac 12+\epsilon$, 
\end{mydef12}

\section{Cross-bred of the Ingham formula and next equivalent of the Fermat-Wiles theorem} 

\subsection{} 

Now we choose the classical Ingham formula, see \cite{2}, 
\be \label{5.1} 
\int_1^T\left|\zf\right|^4{\rm d}t=\frac{1}{2\pi^2}T\ln^4T+\mcal{O}(T\ln^3T),\ T\to\infty
\ee 
as the basic formula for this section. For our purpose we use the following form of the Ingham formula 
\be \label{5.2} 
\int_1^T\left|\zf\right|^4{\rm d}t=\frac{1}{2\pi^2}T\ln^4T\left\{1+\mcal{O}\left(\frac{1}{\ln T}\right)\right\}.
\ee 
Here we give some cross-breeding of the formula (\ref{5.2}) with the following one\footnote{See (\ref{3.11}) for $r=1$.} 
\be \label{5.3} 
\frac{\int_{T}^{\overset{1}{T}}\left|\zf\right|^2{\rm d}t}{\int_{T}^{\overset{1}{T}}|\zeta(\sigma+it)|^2{\rm d}t}=\frac{\ln T}{\zeta(2\sigma)}\left\{1+\mcal{O}\left(\frac{1}{\ln T}\right)\right\}. 
\ee 
Namely, it is true that 
\bdis 
\begin{split}
& \int_1^T\left|\zf\right|^4{\rm d}t= \\ 
& \frac{\zeta^4(2\sigma)}{2\pi^2}T\times 
\left\{
\frac{\int_{T}^{\overset{1}{T}}\left|\zf\right|^2{\rm d}t}{\int_{T}^{\overset{1}{T}}|\zeta(\sigma+it)|^2{\rm d}t}
\right\}^4\times \left\{1+\mcal{O}\left(\frac{1}{\ln T}\right)\right\}
\end{split}
\edis 
or 
\be \label{5.4} 
\begin{split}
& \int_1^T\left|\zf\right|^4{\rm d}t\times \left\{\int_{T}^{\overset{1}{T}}|\zeta(\sigma+it)|^2{\rm d}t\right\}^4\times \\ 
& \left\{\int_{T}^{\overset{1}{T}}\left|\zf\right|^2{\rm d}t\right\}^{-4}=\frac{\zeta^4(2\sigma)}{2\pi^2}T\times \left\{1+\mcal{O}\left(\frac{1}{\ln T}\right)\right\}. 
\end{split}
\ee 

\subsection{} 

Next, we obtain by the substitution 
\be \label{5.5} 
T=\frac{2\pi^2}{\zeta^4(2\sigma)}x\tau,\ \{T\to+\infty\} \ \Leftrightarrow \ \{\tau\to+\infty\}
\ee 
in (\ref{5.4}) the following functional (as cross-bred of the Ingham formula in the preceding sense). 

\begin{mydef54}
\be \label{5.6} 
\begin{split}
& \lim_{\tau\to\infty}\frac{1}{\tau}\int_1^{\frac{2\pi^2}{\zeta^4(2\sigma)}x\tau}\left|\zf\right|^4{\rm d}t\times \\ 
& \left\{\int_{\frac{2\pi^2}{\zeta^4(2\sigma)}x\tau}^{[\frac{2\pi^2}{\zeta^4(2\sigma)}x\tau]^1}|\zeta(\sigma+it)|^2{\rm d}t\right\}^4\times \\ 
& \left\{\int_{\frac{2\pi^2}{\zeta^4(2\sigma)}x\tau}^{[\frac{2\pi^2}{\zeta^4(2\sigma)}x\tau]^1}\left|\zf\right|^2{\rm d}t\right\}^{-4}=x
\end{split}
\ee 
for every fixed 
\bdis 
x>0,\ \sigma\geq \frac 12+\epsilon.
\edis 
\end{mydef54} 

Now, in the special case 
\bdis 
x\to \FR, 
\edis 
we obtain from (\ref{5.6}) the following. 
\begin{mydef55}
\be\label{5.7} 
\begin{split}
& \lim_{\tau\to\infty}\frac{1}{\tau}\int_1^{\frac{2\pi^2}{\zeta^4(2\sigma)}\FR\tau}\left|\zf\right|^4{\rm d}t\times \\ 
& \left\{\int_{\frac{2\pi^2}{\zeta^4(2\sigma)}\FR\tau}^{[\frac{2\pi^2}{\zeta^4(2\sigma)}\FR\tau]^1}|\zeta(\sigma+it)|^2{\rm d}t\right\}^4\times \\ 
& \left\{\int_{\frac{2\pi^2}{\zeta^4(2\sigma)}\FR\tau}^{[\frac{2\pi^2}{\zeta^4(2\sigma)}\FR\tau]^1}\left|\zf\right|^2{\rm d}t\right\}^{-4}=\FR
\end{split}
\ee 
for every fixed Fermat's rational and every fixed $\sigma\geq\frac 12+\epsilon$. 
\end{mydef55} 

Consequently, we have the following result. 

\begin{mydef13}
The $\zeta$-condition 
\be \label{5.8} 
\begin{split}
	& \lim_{\tau\to\infty}\frac{1}{\tau}\int_1^{\frac{2\pi^2}{\zeta^4(2\sigma)}\FR\tau}\left|\zf\right|^4{\rm d}t\times \\ 
	& \left\{\int_{\frac{2\pi^2}{\zeta^4(2\sigma)}\FR\tau}^{[\frac{2\pi^2}{\zeta^4(2\sigma)}\FR\tau]^1}|\zeta(\sigma+it)|^2{\rm d}t\right\}^4\times \\ 
	& \left\{\int_{\frac{2\pi^2}{\zeta^4(2\sigma)}\FR\tau}^{[\frac{2\pi^2}{\zeta^4(2\sigma)}\FR\tau]^1}\left|\zf\right|^2{\rm d}t\right\}^{-4}\not= 1
\end{split}
\ee 
on the class of all Fermat's rationals represents the next $\zeta$-equivalent of the Fermat-Wiles theorem for every fixed $\sigma\geq\frac 12+\epsilon$. 
\end{mydef13}

\section{Cross-bred of the complete Hardy-Littlewood-Ingham formula and non-local interaction of corresponding energies} 

\subsection{} 

In this section we use the following couple of formulas: 
\begin{itemize}
	\item[(a)] Hardy-Littlewood formula\footnote{See (\ref{3.1}).} 
	\be \label{6.1} 
	\zeta(2\sigma)T=\{1+\mcal{O}(T^{-\epsilon}\ln T)\}\int_1^T|\zeta(\sigma+it)|^4{\rm d}t,\ \sigma\geq\frac 12+\epsilon, 
	\ee  
	\item[(b)] our formula (\ref{3.11}), $r=1$ 
	\be \label{6.2} 
	\frac{\ln T}{\zeta(2\sigma)}=\left\{1+\mcal{O}\left(\frac{1}{\ln T}\right)\right\}\times 
	\frac{\int_{T}^{\overset{1}{T}}\left|\zf\right|^2{\rm d}t}{\int_{T}^{\overset{1}{T}}|\zeta(\sigma+it)|^2{\rm d}t},\ \sigma\geq\frac 12+\epsilon
	\ee 
	for triple cross-breeding of the Hardy-Littlewood-Ingham formula (\ref{4.1}). The result is summarized in the next theorem. 
\end{itemize} 

\begin{mydef14}
\be \label{6.3} 
\begin{split}
& \begin{vmatrix}
\int_{T}^{\overset{1}{T}}\left|\zf\right|^2{\rm d}t & 	\int_{T}^{\overset{1}{T}}|\zeta(\sigma+it)|^2{\rm d}t \\ 
\int_{1}^{T}\left|\zf\right|^2{\rm d}t & \int_{1}^{T}|\zeta(\sigma+it)|^2{\rm d}t
\end{vmatrix} \sim \\ 
& \frac{1+\ln2\pi-2c}{\zeta(2\sigma)}\int_{1}^{T}|\zeta(\sigma+it)|^2{\rm d}t\times \int_{T}^{\overset{1}{T}}|\zeta(\sigma+it)|^2{\rm d}t,\ T\to\infty
\end{split}
\ee 
for every fixed $\sigma\geq\frac 12+\epsilon$. 
\end{mydef14} 

\begin{remark}
Bilinear asymptotic formula (\ref{6.3}) expresses a kind of non-local and non-linear interaction between two sets of integrals (i.e. \emph{energies} of corresponding $\zeta$-oscillations), namely 
\be \label{6.4} 
\int_{T}^{\overset{1}{T}}\left|\zf\right|^2{\rm d}t,\ \int_{1}^{T}\left|\zf\right|^2{\rm d}t, 
\ee  
and 
\be \label{6.5} 
\int_{T}^{\overset{1}{T}}|\zeta(\sigma+it)|^2{\rm d}t,\ \int_{1}^{T}|\zeta(\sigma+it)|^2{\rm d}t
\ee 
in the following sense: 
\begin{itemize}
	\item[(a)] 
	\bdis 
	[1,T) \cap (T,\overset{1}{T}]=\emptyset, 
	\edis  
	\item[(b)] 
	\bdis 
	\begin{cases}
	t\in [T,\overset{1}{T}] \ \Rightarrow \ \left\{\frac{1}{2}+it\right\}\cap \{\sigma+it\}=\emptyset, \\ 
	t\in [1,T] \ \Rightarrow \ \left\{\frac{1}{2}+it\right\}\cap \{\sigma+it\}=\emptyset. 
	\end{cases}
	\edis 
\end{itemize}
\end{remark}

I would like to thank Michal Demetrian for his moral support of my study of Jacob's ladders.

\end{document}